\documentclass[12pt]{article}

\input tcilatex

\input xy

\xyoption{all}

\usepackage{amsfonts,epsf}
\usepackage{color}

\newcommand{\psdiag}[3]{\hspace{1mm}\raisebox{-#1mm}{\epsfysize#2mm

\epsffile{#3.eps}}\hspace{1mm}}

\begin{document}

\author{Rui Pedro Carpentier\\
\\ {\small rcarpent@math.ist.utl.pt} \\
{\small\it Departamento de Matem\'{a}tica} \\ {\small\it Centro de
An\'{a}lise Matem\'{a}tica, Geometria e Sistemas Din\^{a}micos}\\
{\small\it Instituto Superior T\'{e}cnico}\\ {\small\it Avenida
Rovisco Pais, 1049-001 Lisboa}\\ {\small\it Portugal}}
\title{On signed diagonal flip sequences\footnote{Supported by Funda\c c\~ao para a Ci\^encia e a Tecnologia, project Quantum Topology, POCI/MAT/60352/2004 and PPCDT/MAT/60352/2004, and project New Geometry and Topology, PTDC/MAT/101503/2008.}}
\date{}
\maketitle

\newpage

\begin{abstract}
Eliahou \cite{2} and Kryuchkov \cite{9} conjectured a proposition that
Gravier and Payan \cite{4} proved to be equivalent to the Four
Color Theorem. It states that any triangulation of a polygon can
be transformed into another triangulation of the same polygon by a
sequence of signed diagonal flips.
It is well known that any pair of polygonal triangulations are
connected by a sequence of (non-signed) diagonal flips. In this
paper we give a sufficient and necessary condition for a diagonal
flip sequence to be a signed diagonal flip sequence. 

\end{abstract}

{\bf keywords:} Four Color Theorem, signed diagonal flip

\newpage

\section{Introduction}

In the study of the Four Color Problem it was observed early on that
the problem could be restricted to graphs that are the 1-skeleton of a 
sphere triangulation, without lost of generality. Due to a result
of Whitney's \cite{10} we can improve the restriction to graphs of this type that
are also Hamiltonian. Thus the Four
Color Theorem is equivalent to the
following statement:

\begin{theorem}\label{xx}
Given two triangulations of a polygon (with no interior vertices) there exists a coloring of
the vertices that is possible for each triangulation.
\end{theorem}

Eliahou \cite{2} and Kryuchkov \cite{9} observed that some moves on
triangulations preserve the coloring (in fact, Kryuchkov worked on
trees instead of polygon triangulations). A proper 4-coloring on the
vertices of a polygon induces a 3-coloring on the edges of the
triangulation such that the boundary of each triangle is colored by
the three colors. This can be done by considering the four colors
as the four elements of the field $\mathbb{F}_4$ of order $4$, and 
then the 3-coloring on the edges of the triangulation is obtained
by coloring each edge with the sum (or difference) of the colors
of its end points. By fixing an order on the three colors we get
a signing on the triangles, $+$ if the colors on the boundary are
ordered in the counterclockwise sense and $-$ if the colors on the
boundary are
ordered in the clockwise sense. In summary, a proper coloring
on the vertices gives a signing for the triangles of a
triangulation of the polygon. If two adjacent triangles have the
same sign then the vertices opposite to the common edge have different
colors and then we can flip the diagonal on the quadrilateral
formed by the two triangles. The sign of the new triangles is the
opposite of the sign of the previous triangles.

$$\psdiag{8}{16}{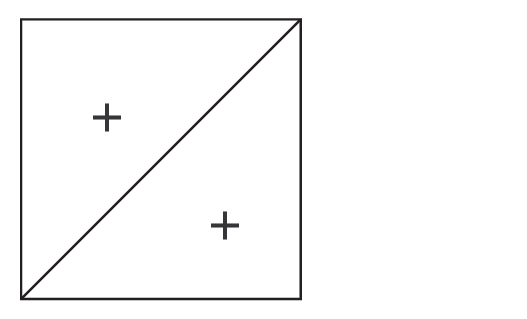} \longleftrightarrow \quad
\psdiag{8}{16}{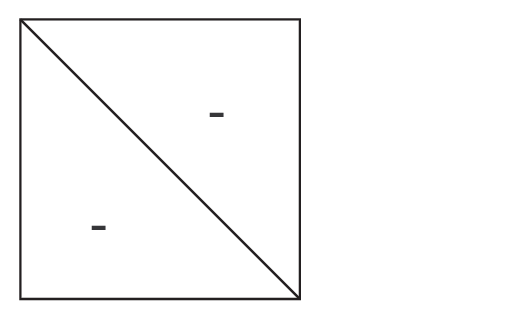}$$

Since this move does not change the coloring on vertices we have
that the existence, for any pair of triangulations of the same
polygon, of signs such that it is possible to transform one to the
other by a sequence of signed diagonal flips implies the Four
Color Theorem, as was observed Eliahou and Kryuchkov. This led
them to conjecture the following:

\begin{conjecture}
Given two triangulations of the same polygon there exist signs for
them such that it is possible to transform one to the other by a
finite sequence of signed diagonal flips.
\end{conjecture}

In \cite{4} Gravier and Payan proved that this
conjecture is, in fact, equivalent to the Four Color Theorem.

It is well known that for any pair of triangulations of a polygon it
is possible to go from one to the other by a finite sequence of
(non-signed) diagonal flips (in the language of binary trees these
are reassociation moves). However, not all sequences of non-signed
diagonal flips can be transformed into a sequence of signed
diagonal flips.

\section{The main result}

We now go on to study when it is possible to transform a sequence
of (non-signed) diagonal flips into a sequence of signed diagonal
flips.

Given a sequence of diagonal flips $\varphi(1)$, $\varphi(2)$, ...
, $\varphi(k)$ from one triangulation of a $n$-polygon to another,
we will construct a graph $G(\varphi)$ in the following way. Given
an enumeration on the vertices of the polygon we represent a
diagonal flip $\varphi(i)$
$$\psdiag{10}{20}{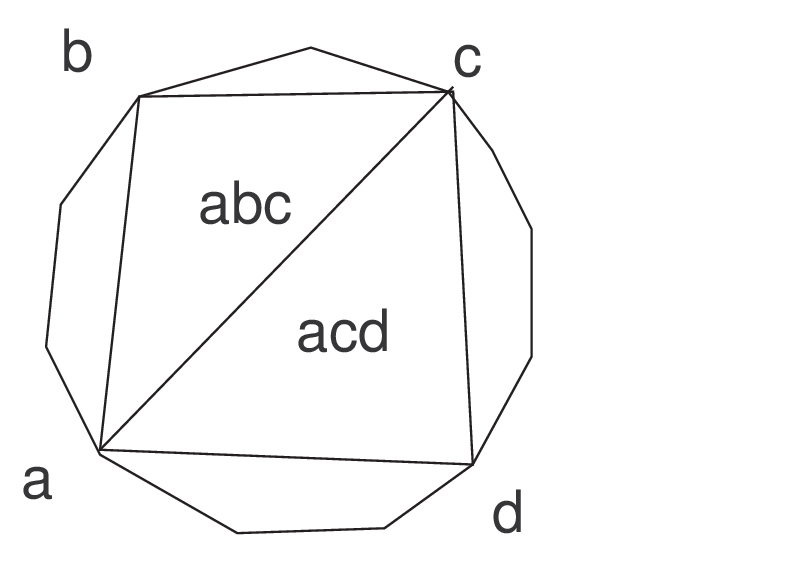}\longrightarrow\psdiag{10}{20}{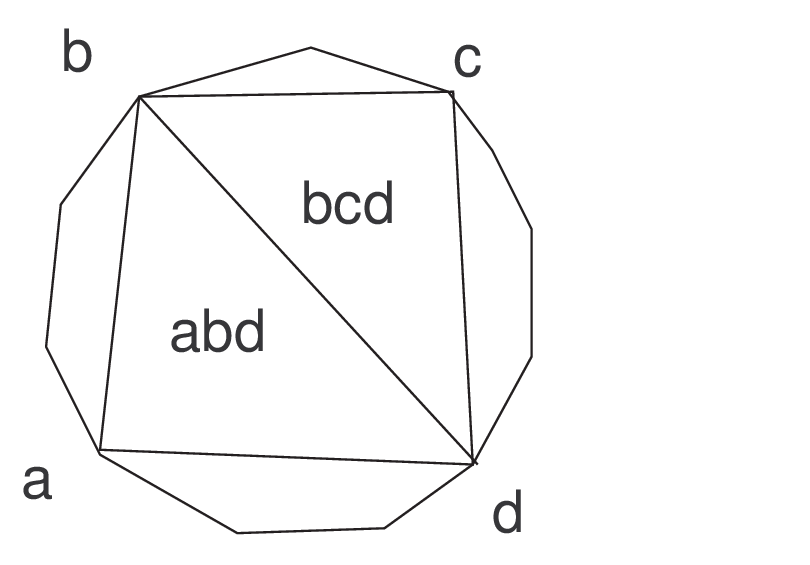}$$ by a
triple $(i,X(i),Y(i))$ where $i$ is the order in the sequence that
$\varphi(i)$ appears, $X(i)=\{abc,acd\}$ is the removed set, i.e. the set of the
triangles in the triangulation that have been removed by the flip
and $Y(i)=\{abd,bcd\}$ is the inserted set, i.e. the set of the triangles in the
triangulation that have been inserted. The diagonal flips will be
the vertices of the graph $G(\varphi)$ and, for $i<j$, the flip
$\varphi(i)=(i,X(i),Y(i))$ is adjacent to the flip
$\varphi(j)=(j,X(j),Y(j))$ if and only if $Y(i)\cap X(j)\not=
\emptyset$ and \begin{equation}Y(i)\cap X(j)\not\subseteq \bigcup_{i<k<j}
X(k).\label{eq:3}
\end{equation} 
Thus, an edge in the graph between the $i^{th}$ and $j^{th}$ vertex occurs 
when one or both triangles in the inserted set $Y(i)$ of the $i^{th}$ flip is 
in the removed set $X(j)$ of the $j^{th}$ flip (there is a "flip interaction"), and 
the same triangle(s) are not involved in an intermediate flip interaction (with flip $k$, for $i<k<j$).
This graph can be constructed in the following way: for $i=1$ to $i=n-1$ we 
take the two triangles in $Y(i)$ (the inserted set of $\varphi(i)$) and join $\varphi(i)$ to 
the first $\varphi(j)$ and the first $\varphi(k)$ (if these exist) for which the respective 
triangles appear in their removed sets $X(j)$ and $X(k)$ ($j$ and $k$ might be the same). 
This algorithm can be carried out in quadratic time on the number of flips (vertices).

{ If we drop the condition (\ref{eq:3}) we get a supergraph of $G(\varphi)$ denoted by $\tilde{G}(\varphi)$.}

{ These graphs give} us a criterion for a sequence of (non-signed)
diagonal flips to be able to be transformed into a sequence of
signed diagonal flips as is shown by the following result.

\begin{theorem}
{The following three statements are equivalent:
\begin{enumerate}
\item[(i)] A sequence $\varphi$ of (non-signed) diagonal flips can be lifted
to a sequence of signed diagonal flips.
\item[(ii)] $\tilde{G}(\varphi)$ is $2$-colorable (i.e. bipartite).
\item[(iii)] $G(\varphi)$ is $2$-colorable.
\end{enumerate}
}
\end{theorem}

\TeXButton{Proof}{\proof}

($(iii)\Rightarrow (i)$) Let $T(k)$ be the set of triangles of the polygonal triangulation 
that precede the flip $\varphi(k)$ and$/$or suceed the flip $\varphi(k-1)$. 
This means that, for each flip $\varphi(k)$, $X(k)\subseteq T(k)$ and $Y(k)\subseteq T(k+1)$. 
Given a $2$-coloring on the graph $G(\varphi)$, for each $k$,
we sign each triangle $t\in T(k)$ by one of the following rules: 

\begin{quote}
R1: search for the last flip $\varphi(i)$ with $i<
k$ such $t\in Y(i)$ and, if it exists, give $t$ the sign $+$ if the color of
$\varphi(i)$ is $1$ and $-$ if the color is $2$;
\end{quote}

\begin{quote}
R2: search for the
first flip $\varphi(j)$ with $j\geq k$ such that $t\in X(j)$ and, if it exists,
give $t$ the sign $-$ if the color of $\varphi(j)$ is $1$ and $+$ if
the color is $2$;
\end{quote}

\begin{quote}
R3: if $t$ doesn't
belong to $X(j)$ for any $j\geq k$ or $Y(i)$ for any $i<k$ we choose {the same sign chosen in $T(k-1)$ 
(if $k=1$ then we choose an arbitrary sign for
$t$)}.
\end{quote}

We need to show that these rules are consistent (i.e. if we apply rules R1 and R2 we obtain the same 
sign for the triangle $t$) and that the signing produced makes the flips become signed flips.

First we observe that if rules R1 and R2 can be both applied {to a triangle $t\in T(k)$ then there 
exist $i$ and $j$ with $i<k\leq j$ such that} \begin{equation} 
t\in Y(i)\cap X(j), \quad t\not\in \bigcup_{i<l<k}
Y(l) \mbox{ and } t\not\in \bigcup_{k\leq l<j}
X(l). \label{eq:1}
\end{equation} 

Also by observing that if a triangle $t$ is not in $T(l)$ nor in $Y(l)$ then it is not in $T(l+1)$ 
{(this is because the new triangles inserted in $T(l+1)$ come from $Y(l)$)}, we can deduce that \begin{equation} 
t\in T(k)\mbox{ and } t\not\in \bigcup_{i<l<k}
Y(l) \quad \Rightarrow t\not\in \bigcup_{i<l<k}
X(l) \label{eq:2}
\end{equation} 
because if $t\in X(l)$ for some $i<l<k$ then $t\not\in T(l+1)$ with $l+1\leq k$ and, if $l+1<k$, we 
can use the fact that $t\not\in Y(l+1)$ for any $i<l+1<k$ to obtain the contradiction $t\not\in T(k)$. 

Thus by putting together (\ref{eq:1}) and (\ref{eq:2}) we have that 
\begin{equation} 
t\in Y(i)\cap X(j)\mbox{ and } t\not\in \bigcup_{i<l<j} X(l)
\end{equation}

$$\mbox{which means that }Y(i)\cap X(j)\not\subseteq \bigcup_{i<l<j} X(l).$$

So $\varphi(i)$ and $\varphi(j)$ have different colors (are adjacent) and rules R1 and R2 induce the same sign on $t$.

Now, let us consider an arbitrary flip $\varphi(i)$ with $X(i)=\{t_1,t_2\}$ and $Y(i)=\{t_3,t_4\}$. 
By rule R2 $t_1$ and $t_2$ receive the same sign ($-$ if the color of
$\varphi(i)$ is $1$ and $+$ if the color is $2$) and by rule R1 $t_3$ and $t_4$ receive the opposite sign 
($+$ if the color of
$\varphi(i)$ is $1$ and $-$ if the color is $2$). Thus $\varphi(i)$ becomes a signed flip.

{We note that if a triangle $t\in T(k)$ is not used in the flip $\phi(k)$ then it keeps the same sign after 
the flip by the rules $R1$, $R2$ and $R3$.}

{($(ii)\Rightarrow (iii)$) This follows immediately from the fact that $G$ is a subgraph of $\tilde{G}$.}

{($(i)\Rightarrow (ii)$) If a sequence of diagonal flips $\varphi$ is a sequence of signed
diagonal flips then this induces a $2$-coloring on the graph
$G(\varphi)$ by coloring each vertex $\varphi(i)$ with the sign of
the triangles on its removed set $Y(i)$. To see that this $2$-coloring is in fact well defined we observe that, 
by the result of Eliahou and Kryuchkov, there exists a $4$-coloring on the vertices of the polygon that remains 
valid and invariant under the sequence of flips, and the signs of the triangles are determined by the colors of 
these vertices. Thus a fixed triangle has the same sign in any triangulation to which it belongs, so if two 
vertices $\varphi(i)$ and $\varphi(j)$ are adjacent ($Y(i)$ and $X(j)$ share a triangle $t$) then they are colored 
with different colors ($\varphi(i)$ is colored with the sign of $t$ and $\varphi(j)$ is colored with the opposite 
sign since $X(j)$ and $Y(j)$ have opposite signs).} \TeXButton{End Proof}{\endproof}

\section{Some examples and remarks}

In the first two examples both graphs $G$ and $\tilde{G}$ are equal.

\begin{example}
Consider the following sequence of flips:

$$\psdiag{10}{20}{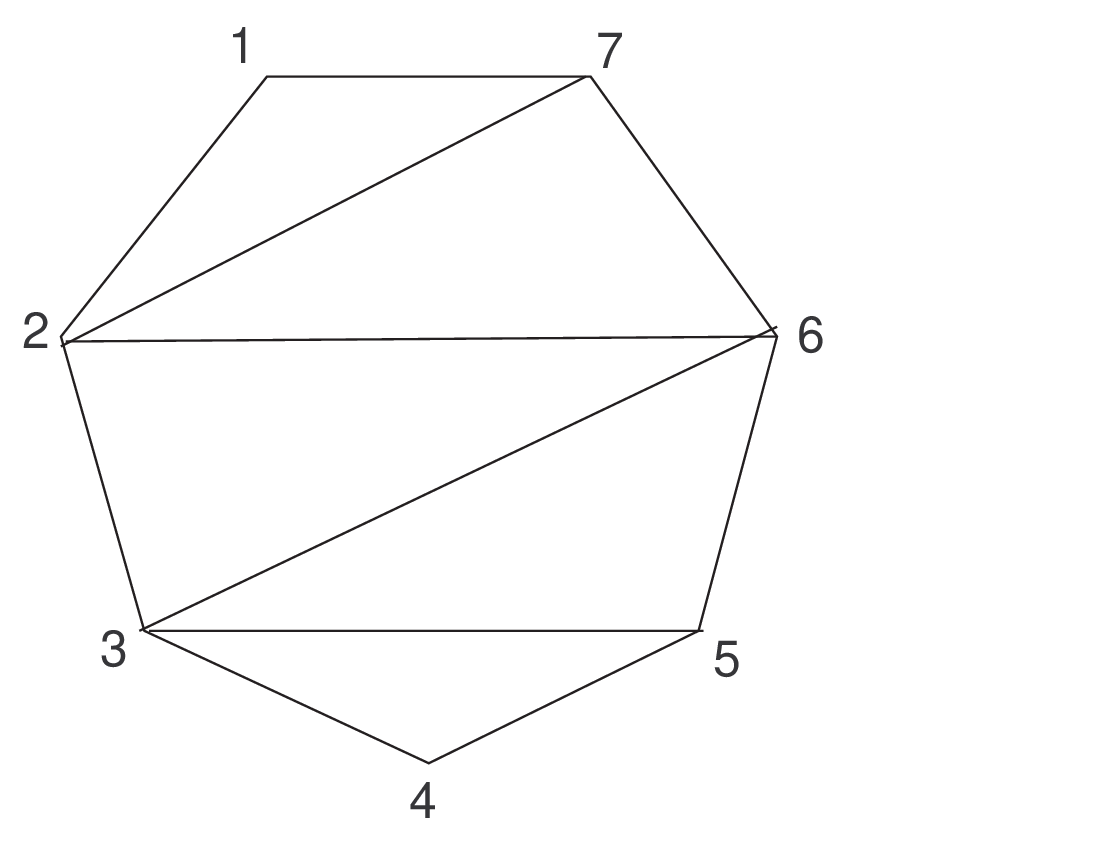}\rightarrow \psdiag{10}{20}{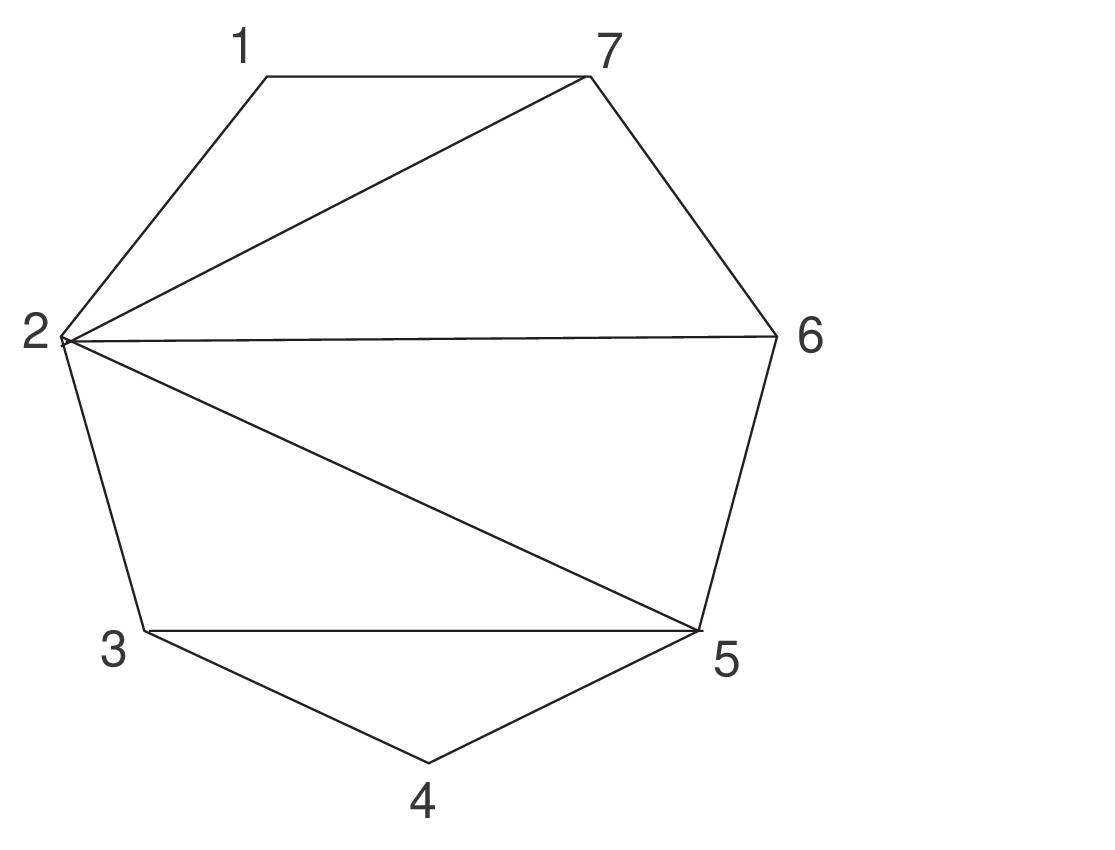}\rightarrow
\psdiag{10}{20}{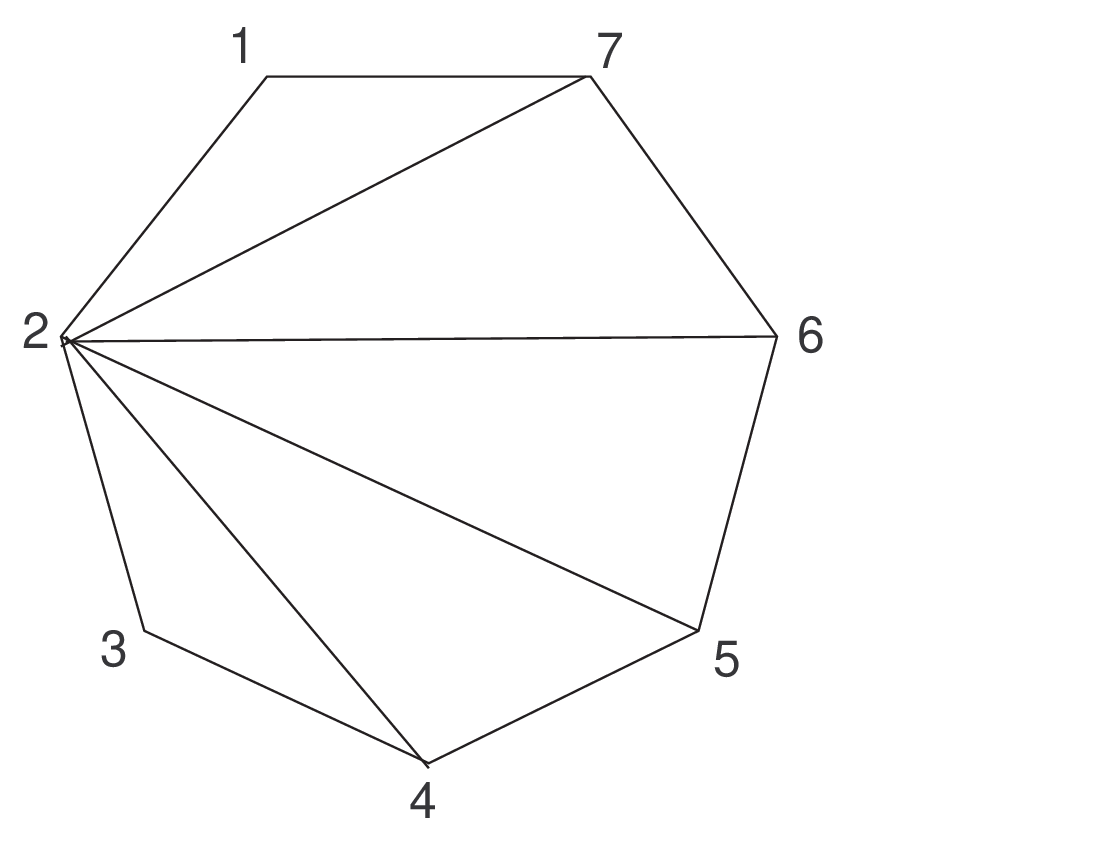}\rightarrow \psdiag{10}{20}{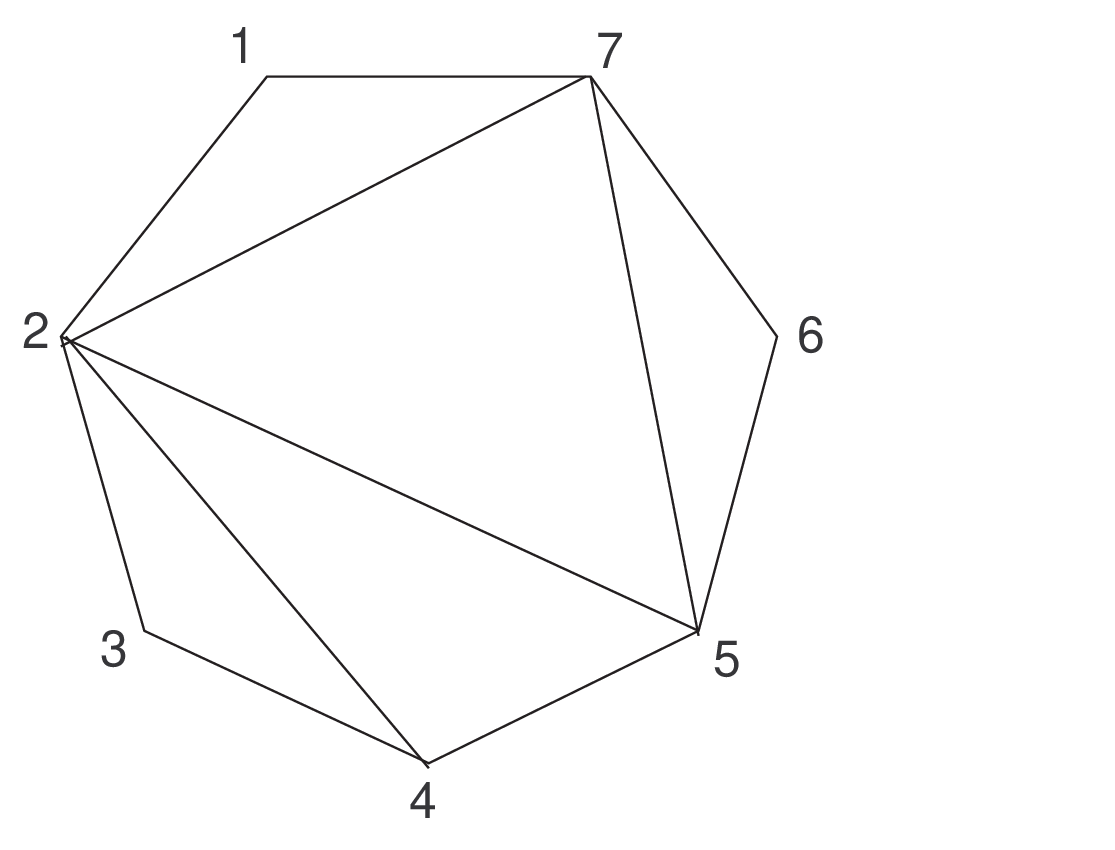}$$ $$
\rightarrow \psdiag{10}{20}{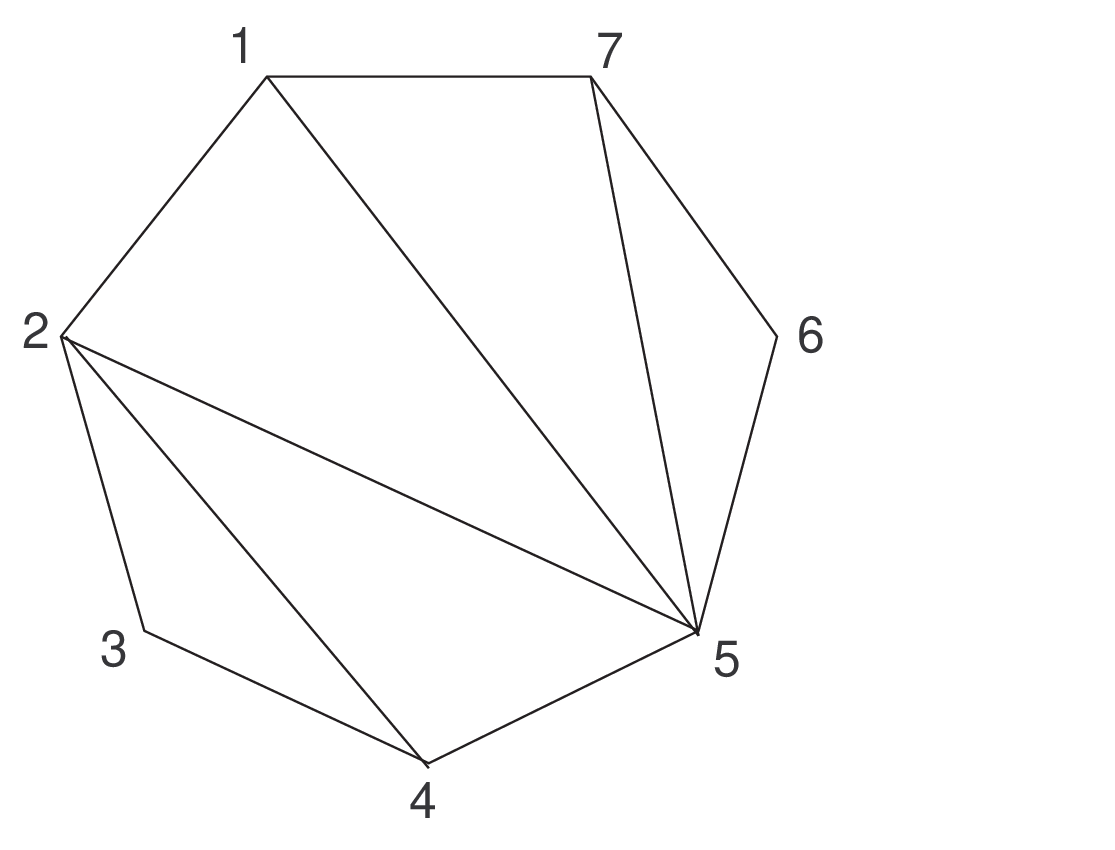}\rightarrow
\psdiag{10}{20}{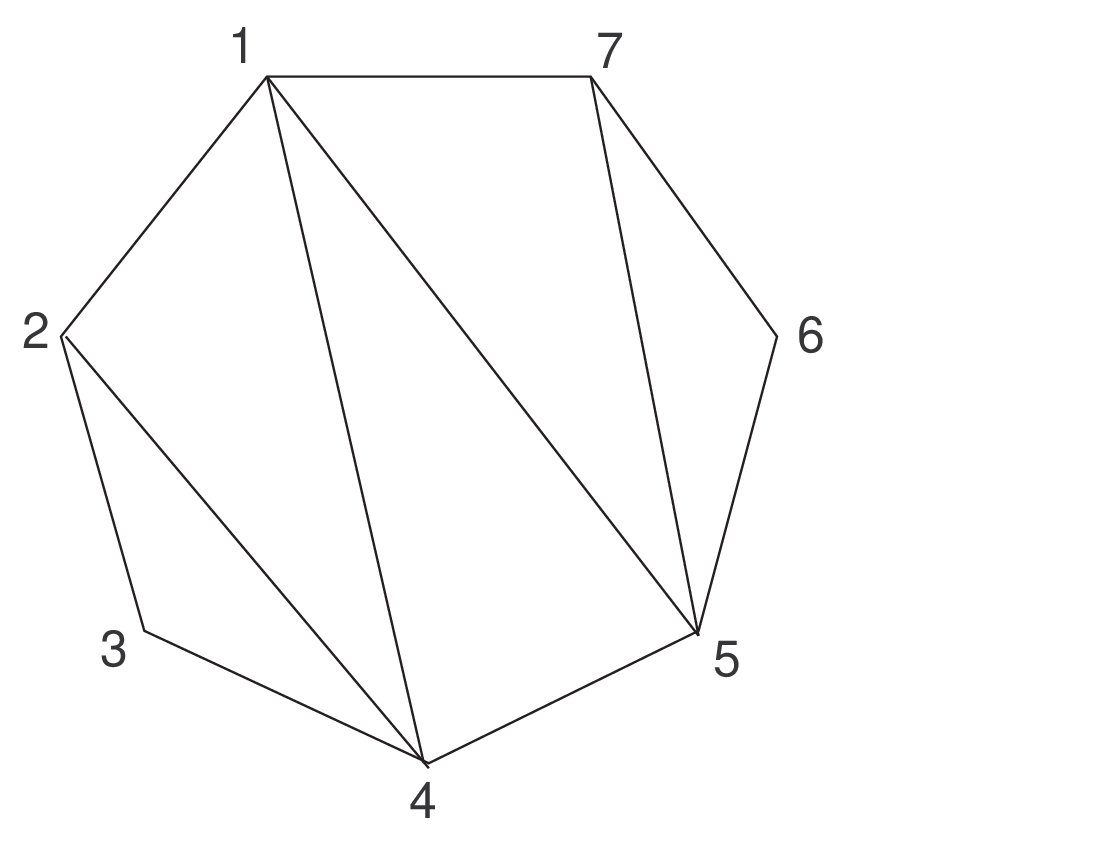}$$

$\varphi(1)=(1,\{236,356\},\{235,256\})$

$\varphi(2)=(2,\{235,345\},\{234,245\})$

$\varphi(3)=(3,\{256,267\},\{257,567\})$

$\varphi(4)=(4,\{127,257\},\{125,157\})$

$\varphi(5)=(5,\{125,245\},\{124,145\})$

This gives a non-$2$-colorable graph:

$$\xymatrix{&\varphi(1) \ar @{-}[dl] \ar @{-}[rr]&
&\varphi(2)\ar@{-}[dr]&\\ \varphi(3)\ar @{-}[drr]& & & &
\varphi(5)\ar @{-}[dll] \\ & & \varphi(4) & & }$$

Thus the sequence of flips is not signable.
\end{example}

\begin{example}

Now we consider another sequence of flips:

$$\psdiag{10}{20}{hep1}\rightarrow \psdiag{10}{20}{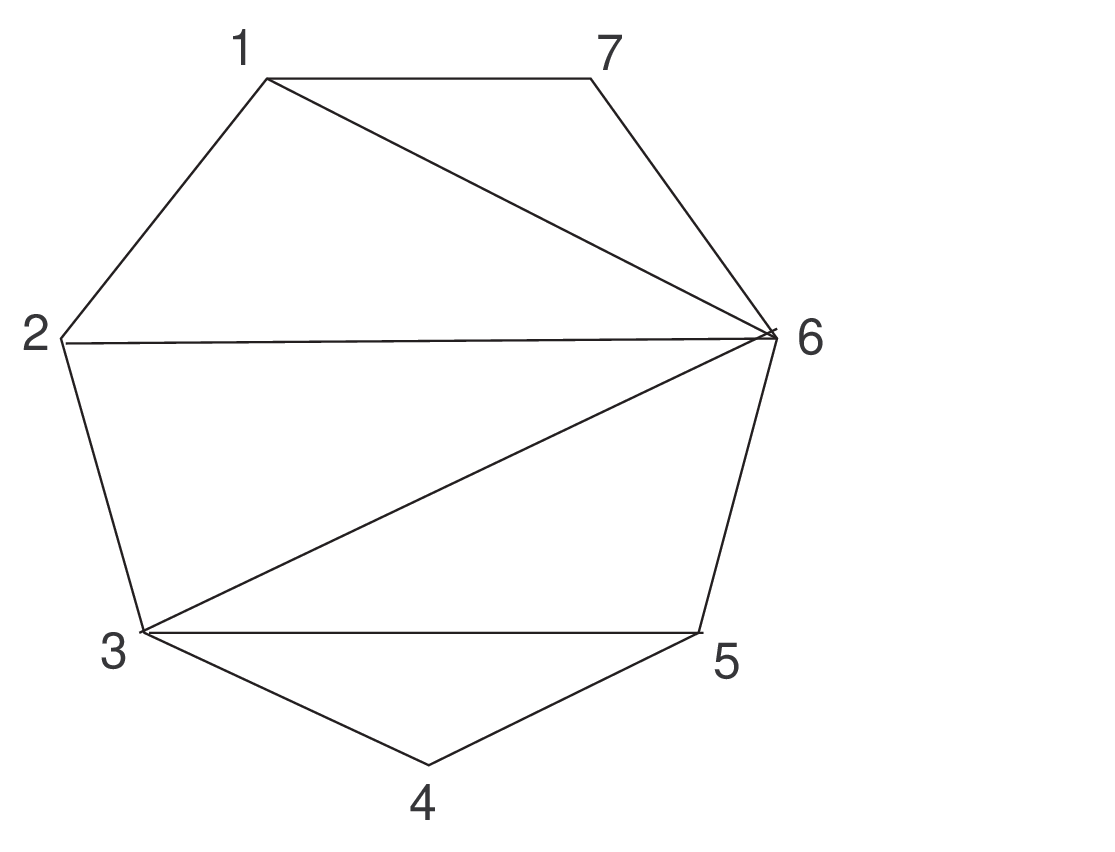}\rightarrow
\psdiag{10}{20}{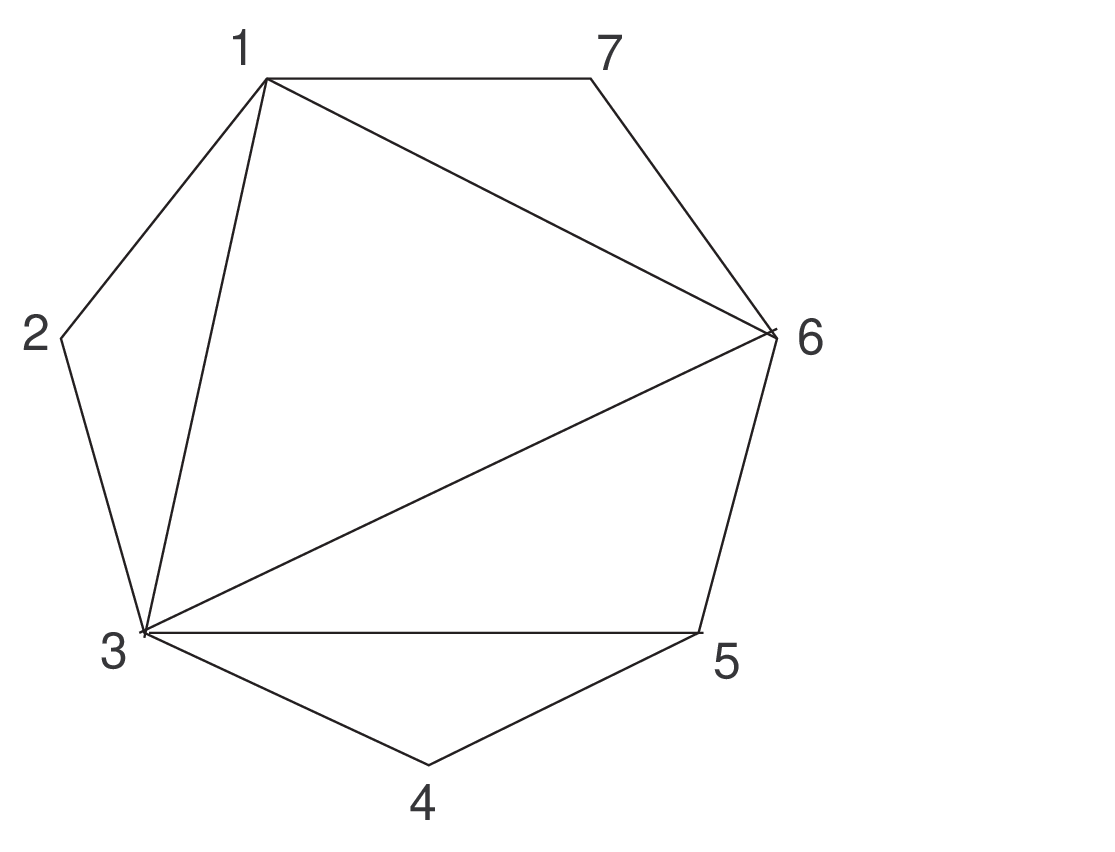}\rightarrow \psdiag{10}{20}{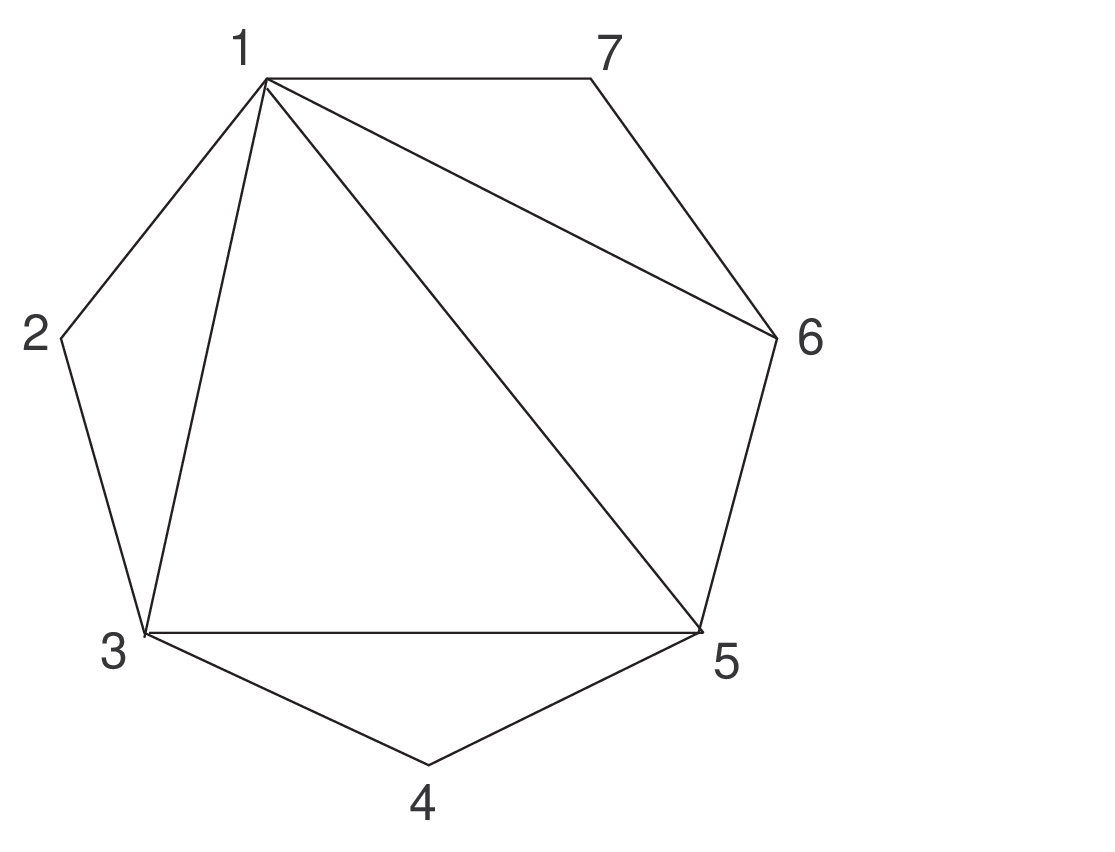}$$ $$
\rightarrow \psdiag{10}{20}{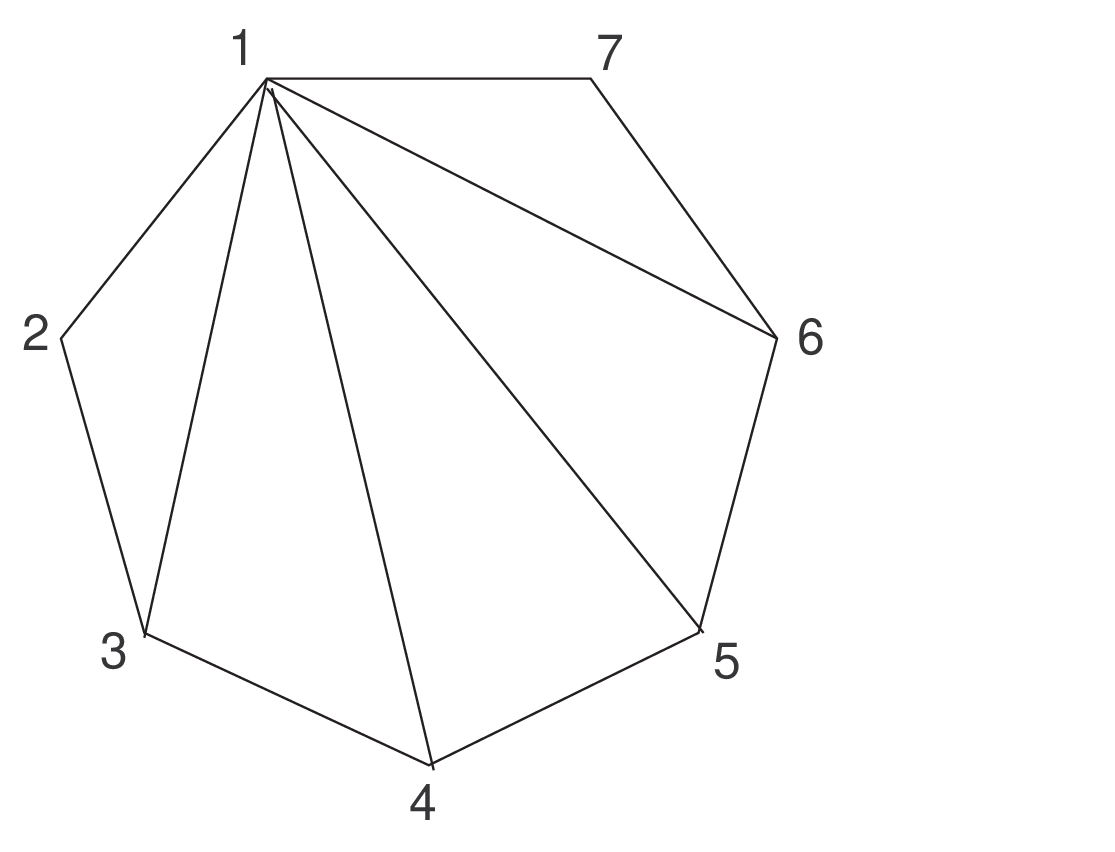}\rightarrow
\psdiag{10}{20}{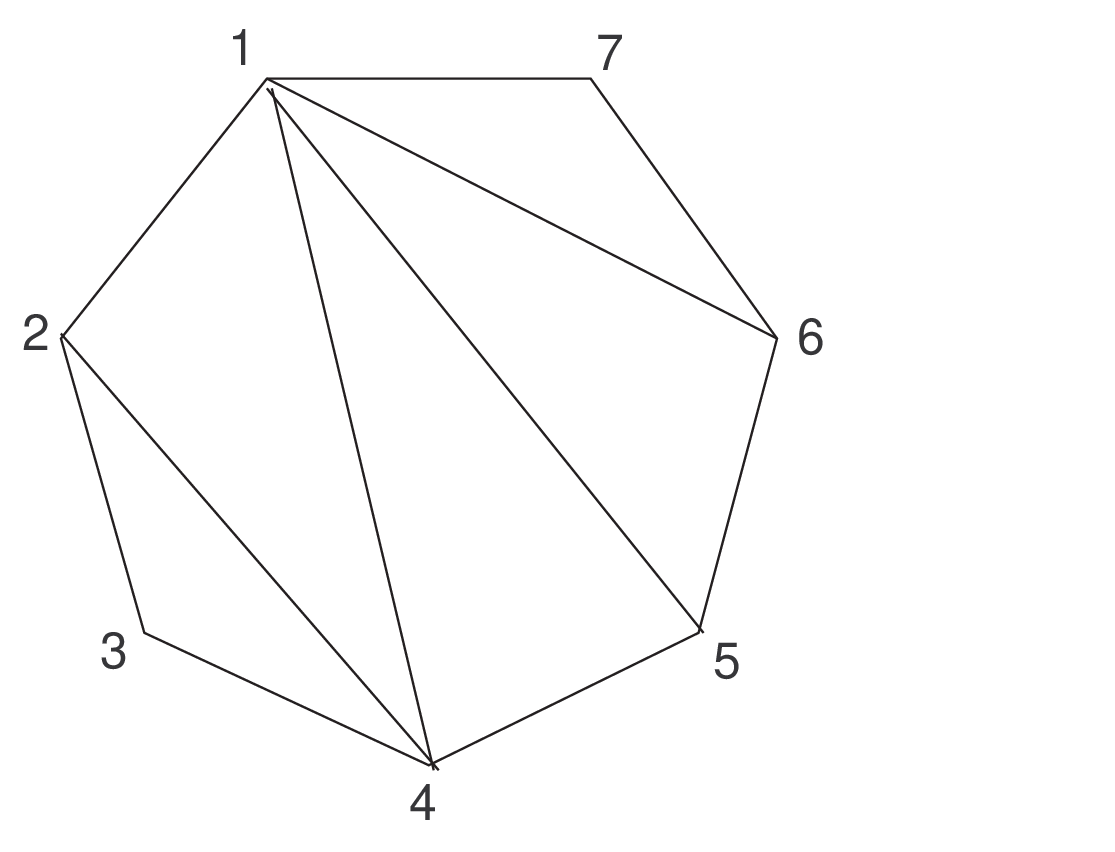}\rightarrow \psdiag{10}{20}{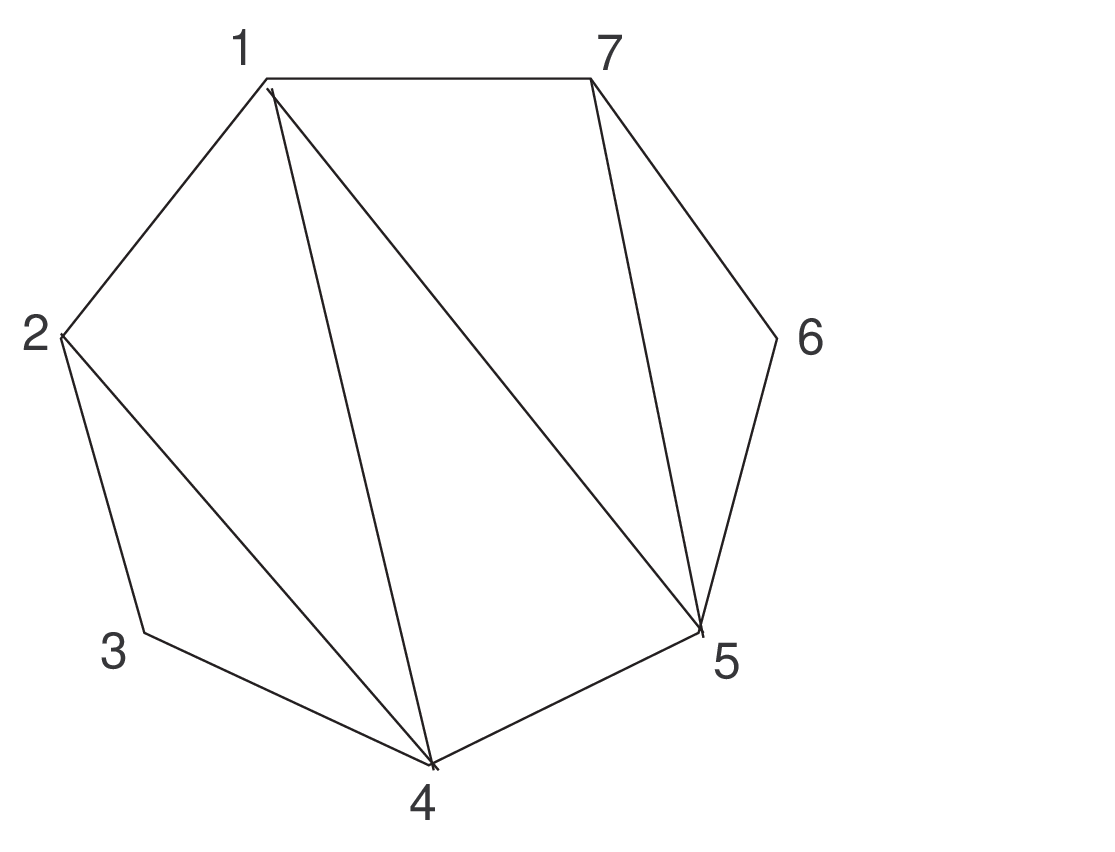}$$

$\varphi(1)=(1,\{127,267\},\{126,167\})$

$\varphi(2)=(2,\{126,236\},\{123,136\})$

$\varphi(3)=(3,\{136,356\},\{135,156\})$

$\varphi(4)=(4,\{135,345\},\{134,145\})$

$\varphi(5)=(5,\{123,134\},\{124,234\})$

$\varphi(6)=(6,\{167,156\},\{157,567\})$

This gives a $2$-colorable graph:

$$\xymatrix{\varphi(2)\ar @{-}[rrr]\ar @{-}[dd]\ar @{-}[dr]& & &
\varphi(3)\ar @{-}[dl]\ar @{-}[dd]
\\&\varphi(5) \ar @{-}[r] \ar @{-}[r]&\varphi(4)&\\
 \varphi(1)\ar @{-}[rrr] & & & \varphi(6)}$$
which gives a sequence of signed flips:

$$\psdiag{10}{20}{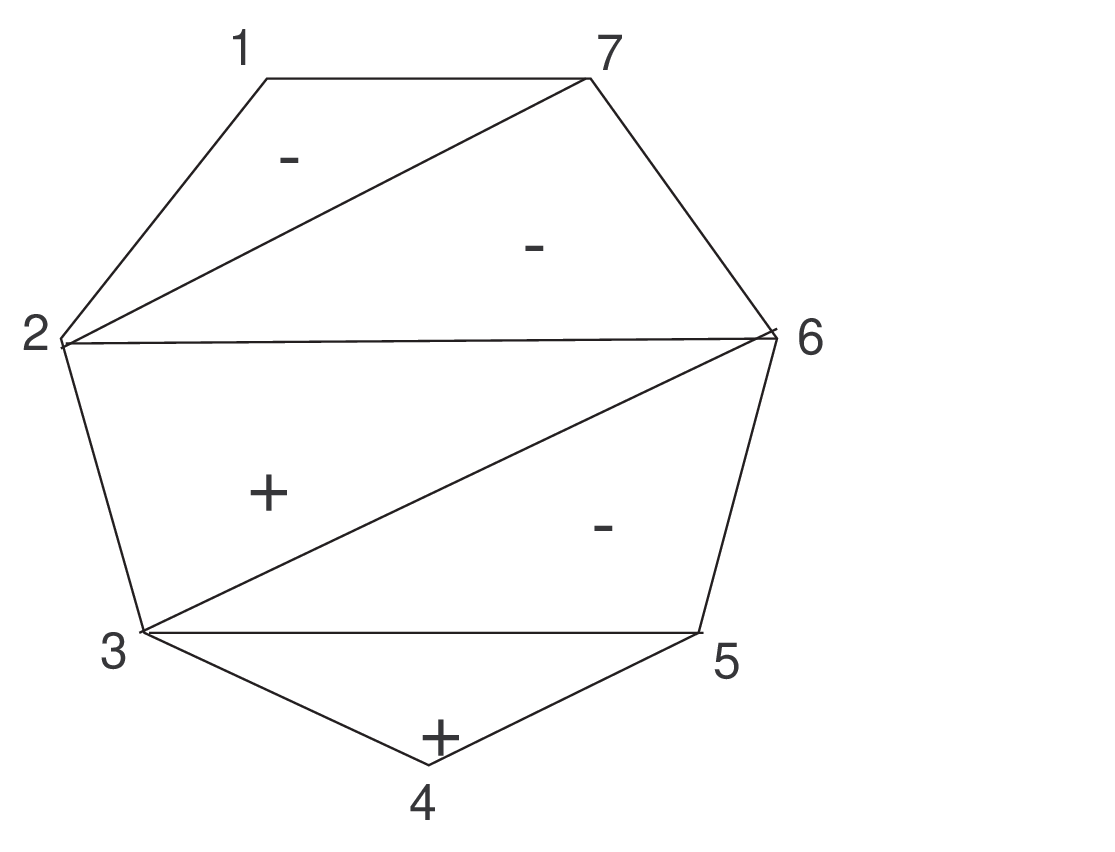}\rightarrow
\psdiag{10}{20}{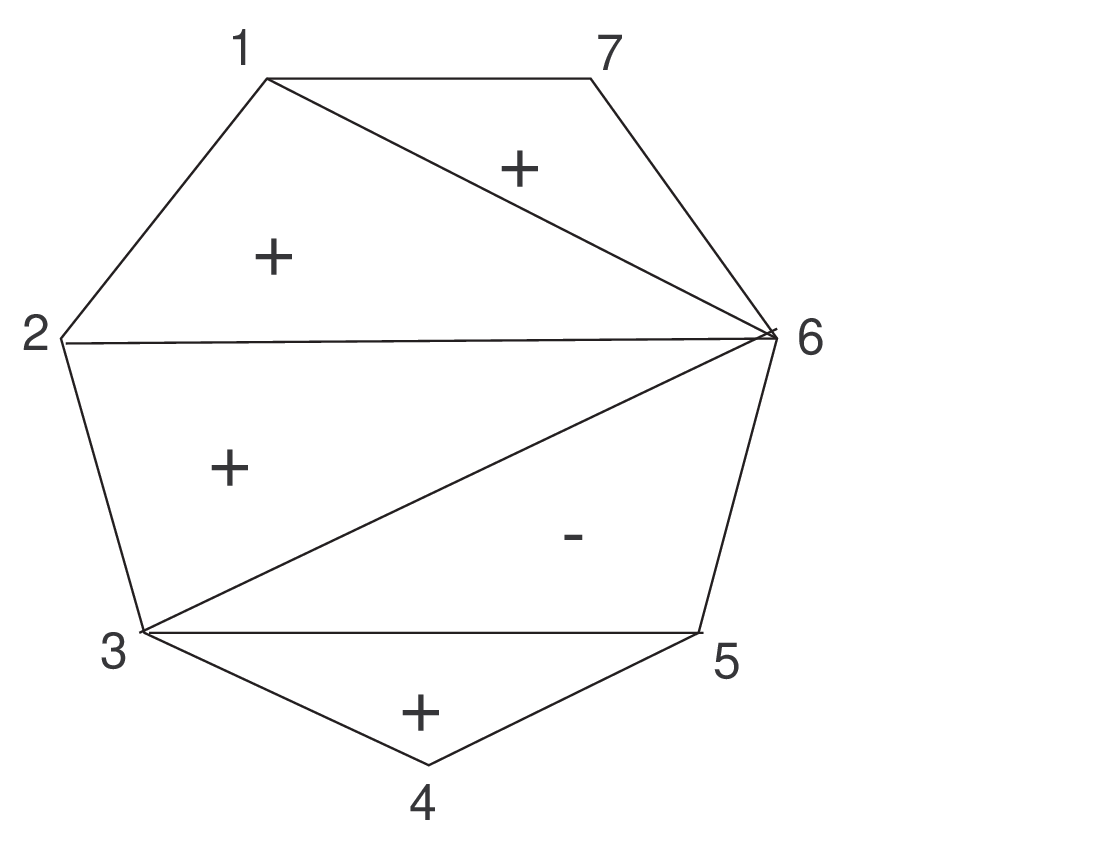}\rightarrow \psdiag{10}{20}{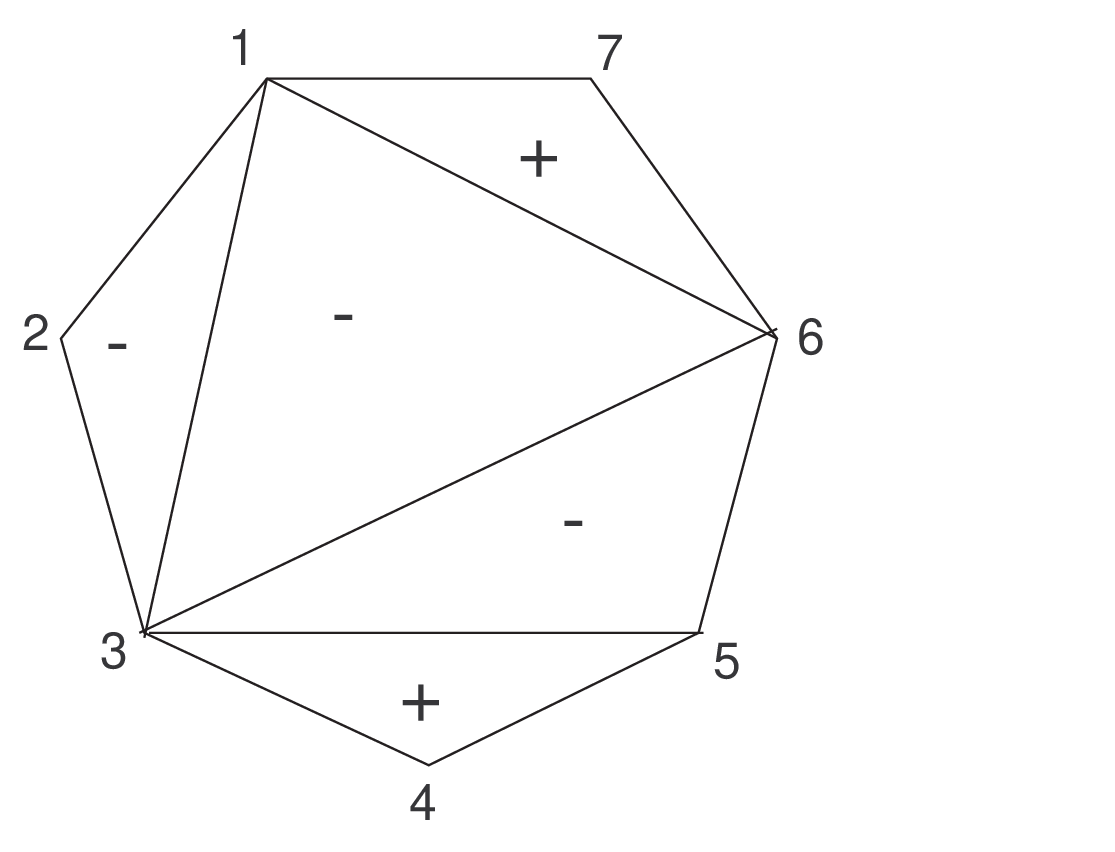}\rightarrow
\psdiag{10}{20}{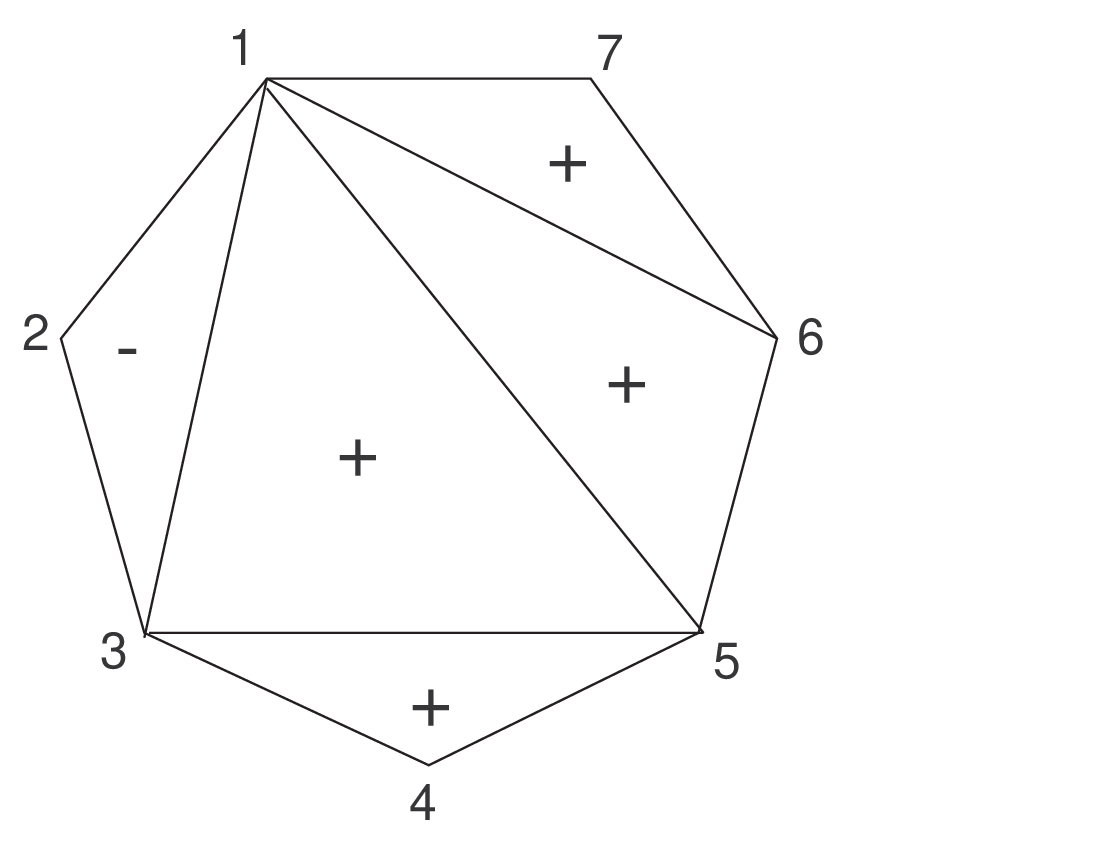}$$ $$\rightarrow
\psdiag{10}{20}{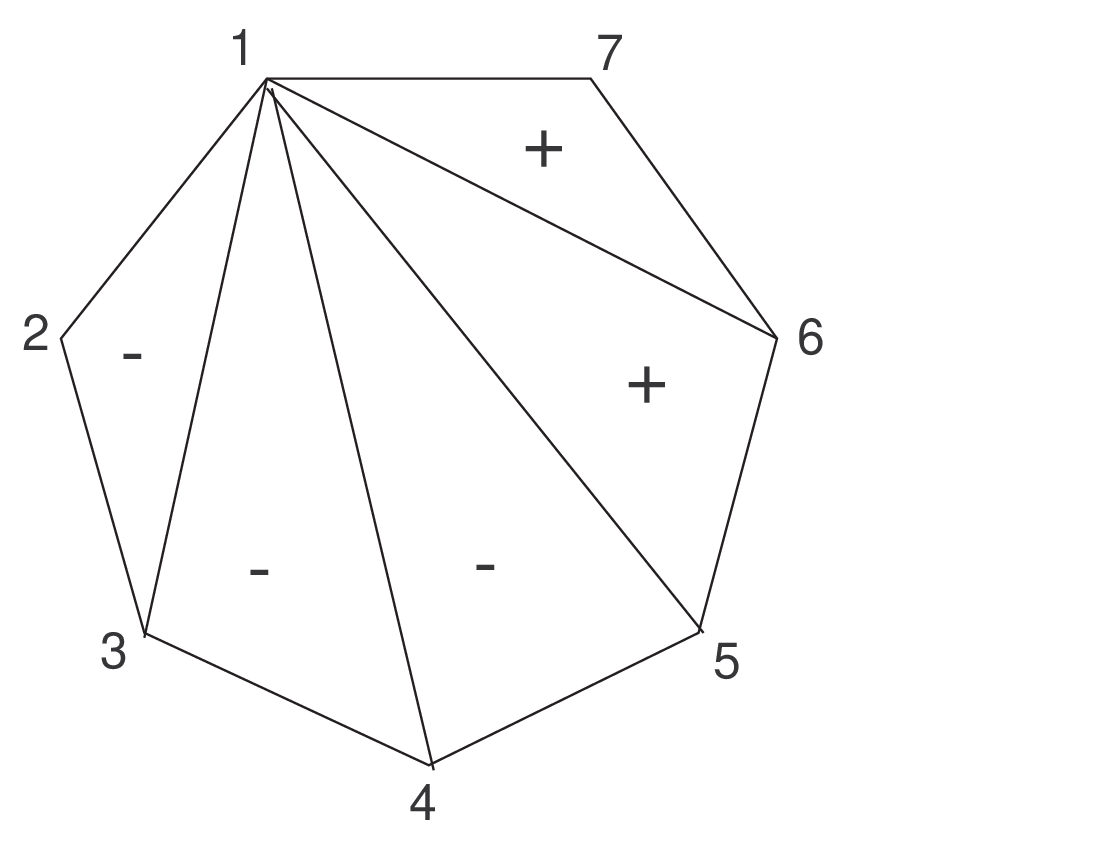}\rightarrow \psdiag{10}{20}{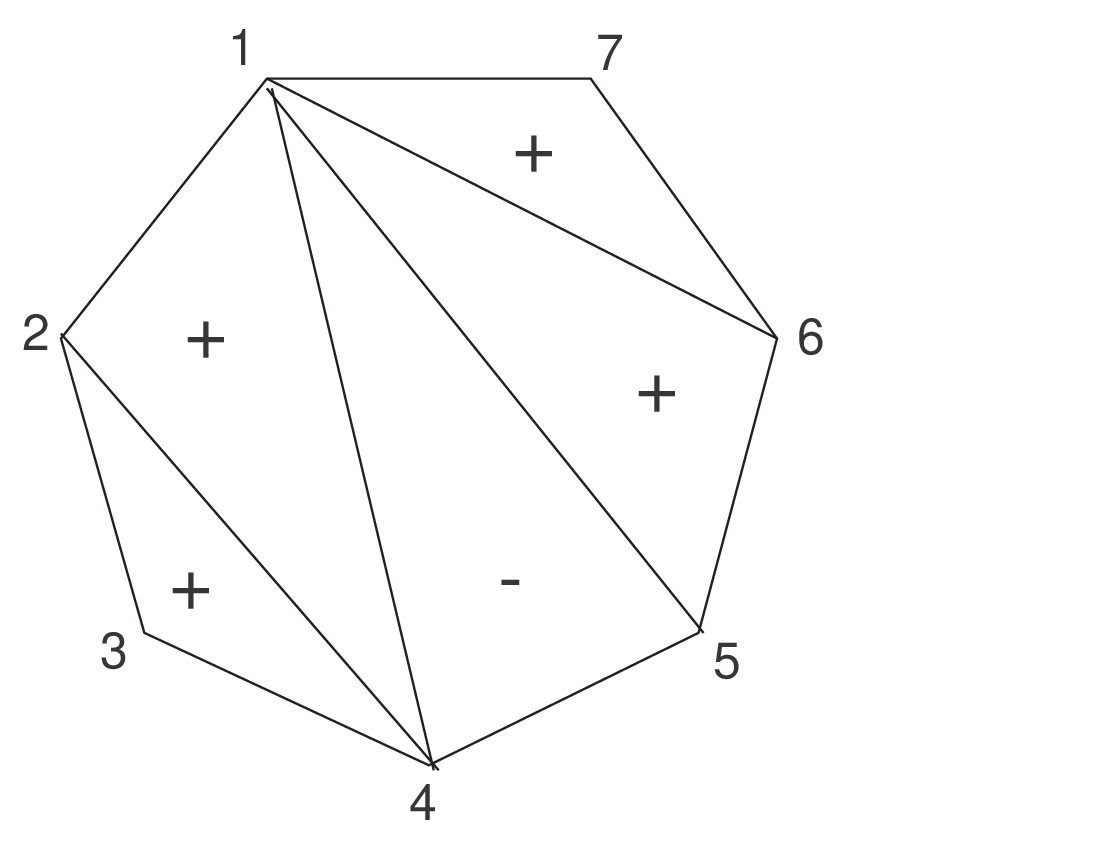}\rightarrow
\psdiag{10}{20}{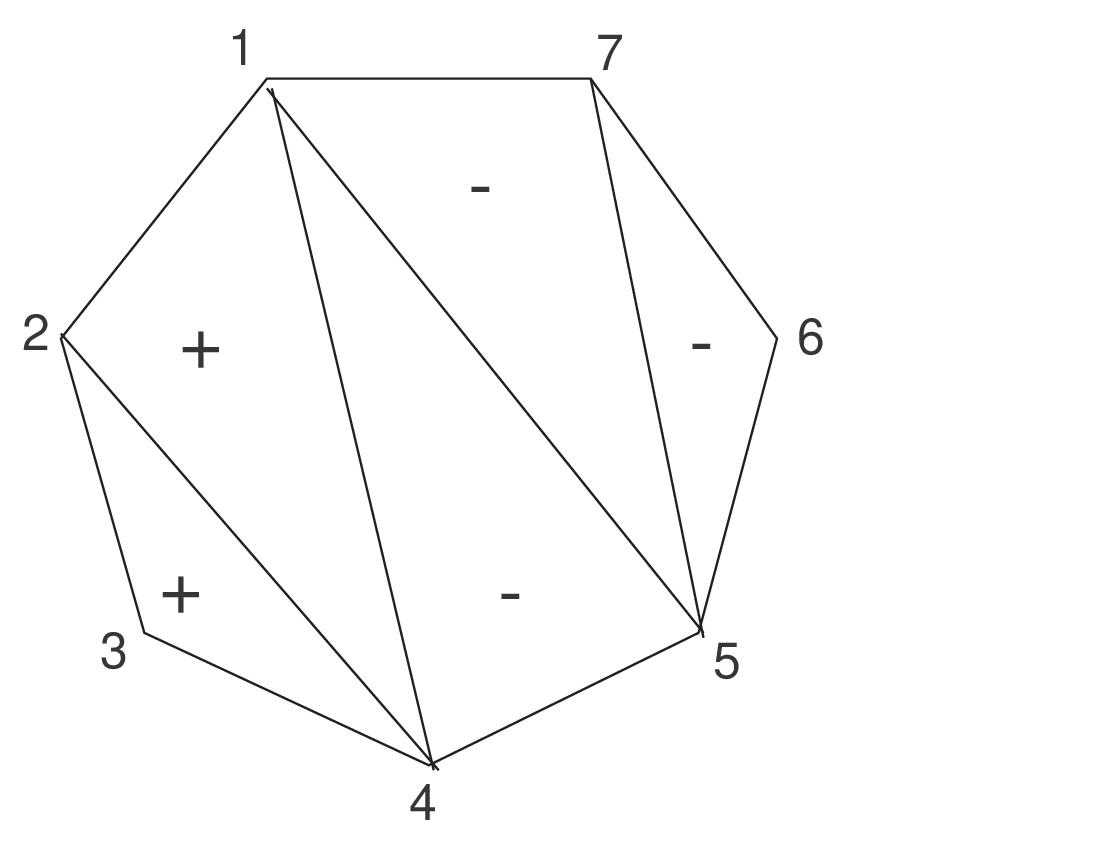}$$
\end{example}

Although the graph $\tilde{G}$ has a simpler definition and coincides with $G$ for 
``small'' sequences of flips, they are not the same as can be seen in the following example.

\begin{example}

Consider the following sequence of flips:

$$\psdiag{10}{20}{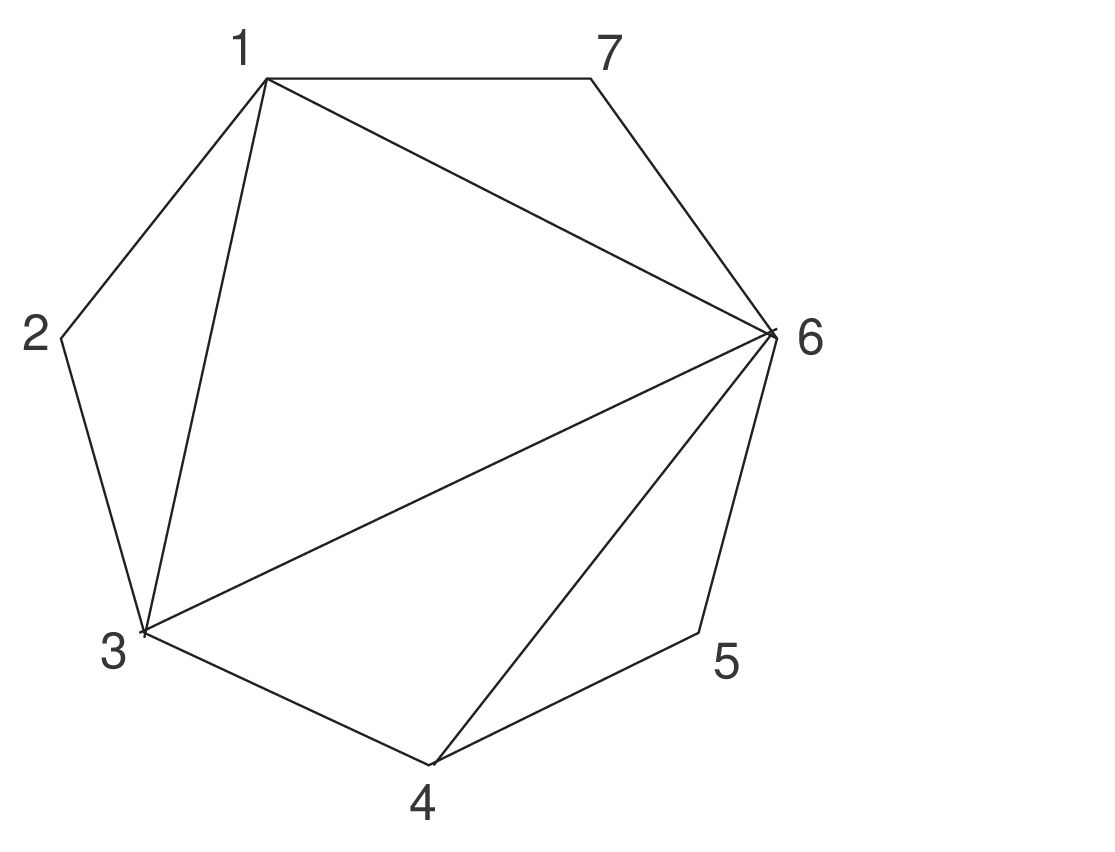}\rightarrow \psdiag{10}{20}{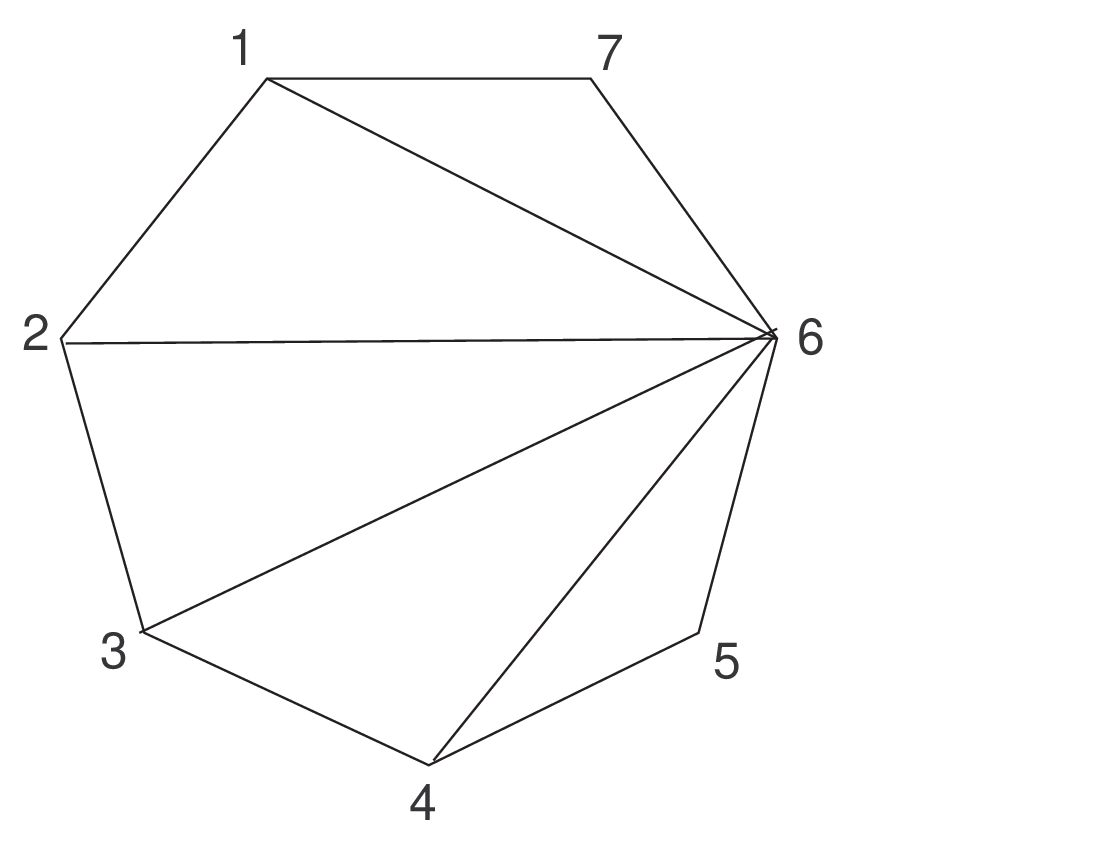}\rightarrow
\psdiag{10}{20}{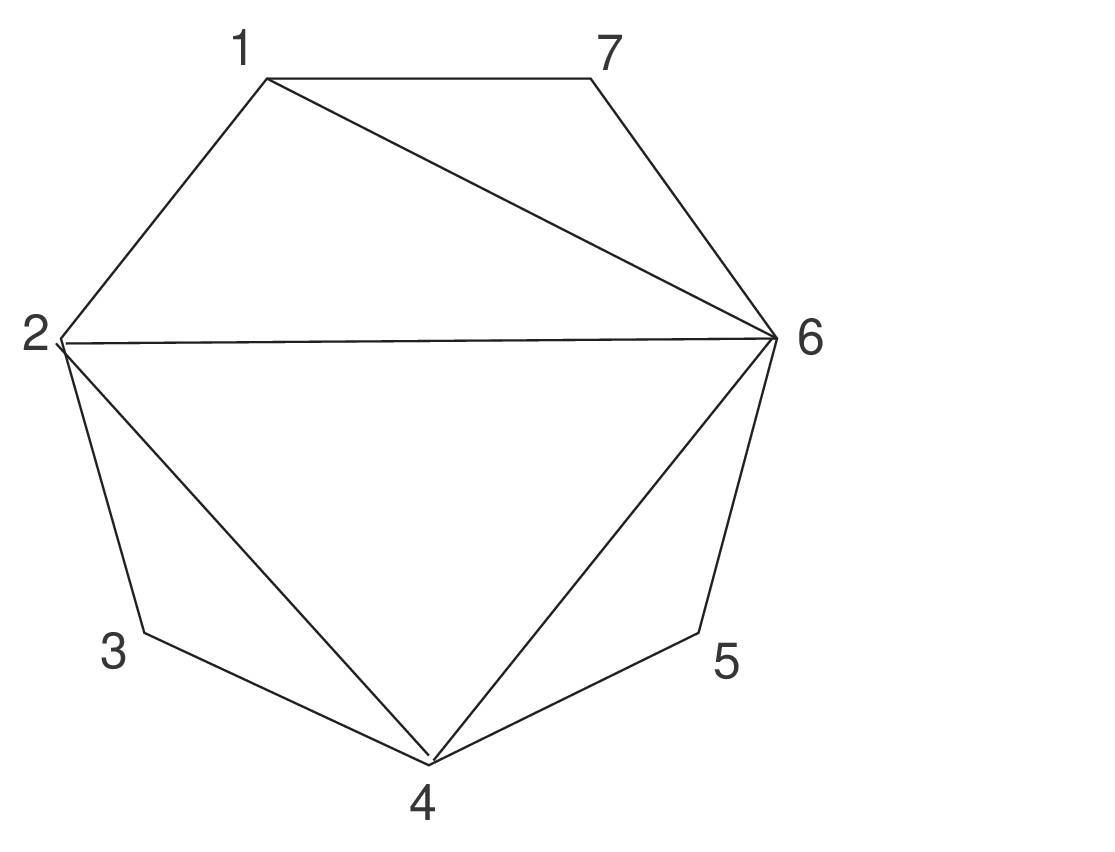}\rightarrow \psdiag{10}{20}{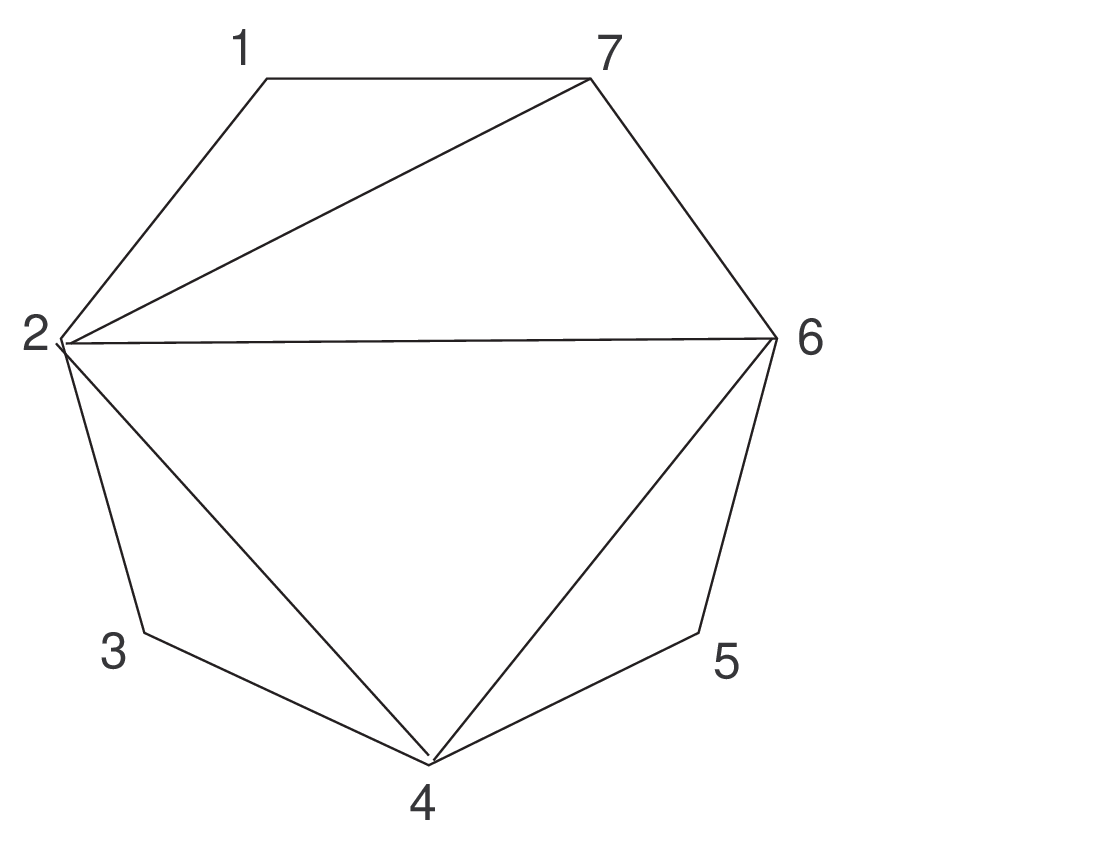}$$ $$
\rightarrow \psdiag{10}{20}{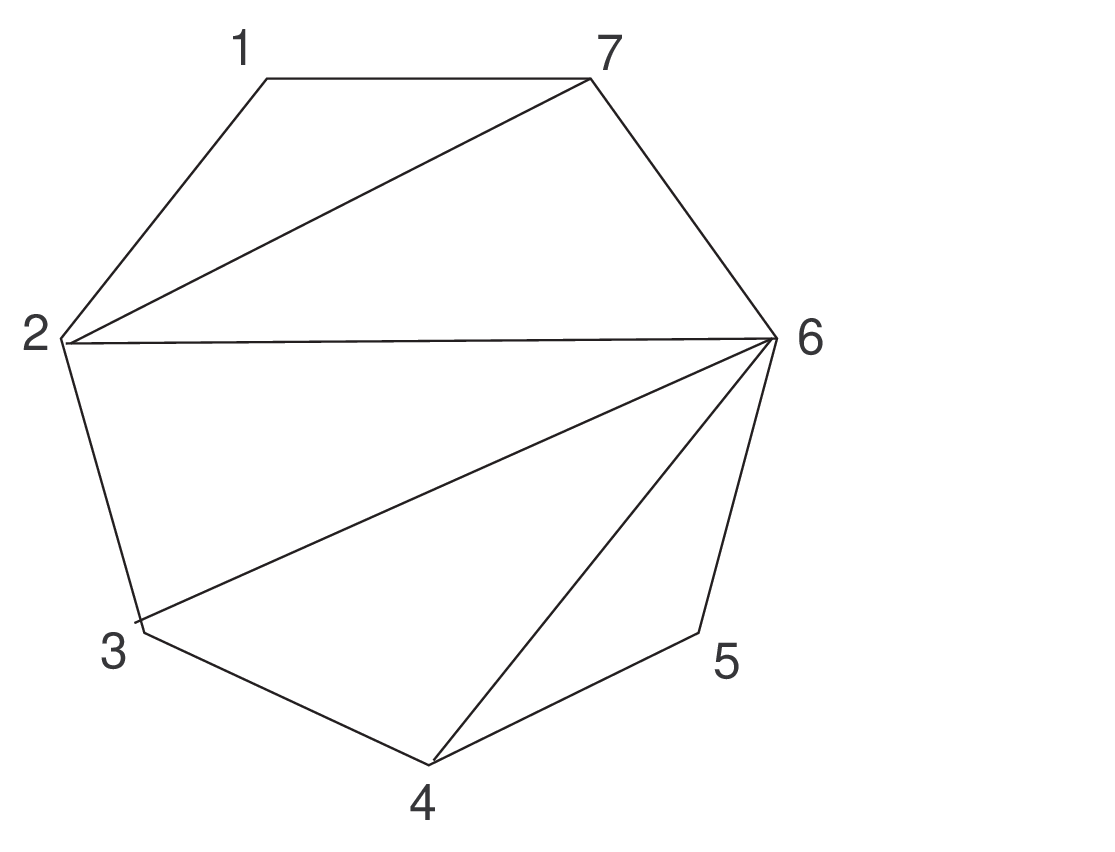}\rightarrow
\psdiag{10}{20}{hep1}\rightarrow
\psdiag{10}{20}{hep2b}\rightarrow
\psdiag{10}{20}{hep3b}$$

$\varphi(1)=(1,\{123,136\},\{126,236\})$

$\varphi(2)=(2,\{236,345\},\{234,246\})$

$\varphi(3)=(3,\{126,167\},\{127,267\})$

$\varphi(4)=(4,\{234,246\},\{236,346\})$

$\varphi(5)=(5,\{346,456\},\{345,356\})$

$\varphi(6)=(6,\{127,267\},\{126,167\})$

$\varphi(7)=(7,\{126,236\},\{123,136\})$

This gives the two $2$-colorable graphs $G$ and $\tilde{G}$ (the dashed line is an edge that belongs to $\tilde{G}$ but not to $G$):

$$\xymatrix{\varphi(1)\ar @{-}[r]\ar @{-}[d]\ar @{.}[drr]& \varphi(3)\ar @{-}[r] & \varphi(6) \ar @{-}[d]\\
\varphi(2)\ar @{-}[r] &\varphi(4) \ar @{-}[r] \ar @{-}[d]&\varphi(7)\\
 & \varphi(5)& }$$

\end{example}

We observe in these examples that the degree of each vertex of $G$ is at most three. In general, we have that this degrees is at most four, which follows from the non-obvious fact that the graph $G$ is the same (up to isomorphism) as the graph obtained after reversing the order of the flip sequence. We note that the same symmetry holds trivially for $\tilde{G}$.

The reason is that, since each vertex $\varphi(i)$ in $G$ is adjacent to at most two vertices $\varphi(j)$ and $\varphi(k)$ with indices greater than $i$ (at most one for each triangle in $Y(i)$) by this symmetry the same vertex is adjacent to at most two vertices $\varphi(j')$ and $\varphi(k')$ with indices smaller than $i$ (at most one for each triangle in $X(i)$), thus the degree of $\varphi(i)$ is at most four. To see this invariance under reversing the order of a flip sequence, we observe that, for $i<j$, if a triangle is in $Y(i)\cap X(j)$ (which means that the triangle is inserted in the flip $\varphi(i)$) but it is not in $\bigcup_{i<k<j}X(k)$ (which means that it is not removed by any flip $\varphi(k)$) then it is in all triangulations of the polygon between the flips $\varphi(i)$ and $\varphi(j)$, thus it cannot be in $\bigcup_{i<k<j}Y(k)$ (we cannot insert a triangle that is already in the triangulation). On the other hand, if a triangle is in $Y(i)\cap X(j)$ (which means that the triangle is removed in the flip $\varphi(j)$) but it is not in $\bigcup_{i<k<j}Y(k)$ (which means that it is not inserted after the flip $\varphi(i)$) then it cannot be in $\bigcup_{i<k<j}X(k)$ (if it is removed by some flip $\varphi(k)$ then it must be inserted by some flip $\varphi(k')$ with $k<k'<j$). Therefore $Y(i)\cap X(j)$ is not a subset of $\bigcup_{i<k<j}X(k)$ if and only if it is not a subset of $\bigcup_{i<k<j}Y(k)$. This implies that the graph $G$ of the inverted flip sequence is the same as the graph $G$ of the original sequence, since the removed sets of the inverted sequence are the inserted sets of the original sequence and vice-versa.

\vspace{1cm}

{\bf Acknowledgments} - I wish to thank Roger Picken for his useful suggestions and comments. {I also wish to thank the anonymous referee whose suggestions and comments helped to improve this manuscript.}

\end{document}